\documentclass[12pt]{extarticle}

\usepackage[english]{babel}
\usepackage{amssymb,amsthm,amsmath,mathtools}
\usepackage{subcaption}
\usepackage{tikz}
\usepackage{graphicx}
\usetikzlibrary{calc}
\usetikzlibrary{decorations.pathreplacing}
\tikzstyle{snode}=[circle ,draw=black,fill=red,thick, inner sep=0pt ,minimum size=1.4mm]
\tikzstyle{bnode}=[circle ,draw=black,fill=black,thick, inner sep=0pt ,minimum size=1.4mm]
\usepackage{textcomp}
\usepackage{ifthen}
\usepackage[normalem]{ulem}

\usepackage{enumitem}
\setlist[enumerate]{left=0pt,label=(\roman*),
ref=(\roman*),font=\normalfont,topsep=-1ex,parsep=-.3ex,partopsep=0pt}

\usepackage{hyperref}
\hypersetup{
    colorlinks=true,
    pdftitle={Lower bounds on the indepndence number of a graph in terms of degrees},
    pdfauthor={J. Harant, I.Schiermeyer},
    bookmarks=true,
    pdfpagemode=UseOutlines,
}
\usepackage[capitalise]{cleveref}

\usepackage[left=2.5cm,right=2.5cm,top=3cm,bottom=3cm]{geometry} 
\newtheorem{theorem}{Theorem}

\newtheorem{lemma}{Lemma}

\title{{\bf Lower bounds on the independence number of a graph in terms of degrees}}

\author{Jochen Harant$^1$ and Ingo Schiermeyer$^{2,3}$  \\~\\
\small $^1$  Technical University of  Ilmenau, 98693 Ilmenau, 
\small Germany\\
\small $^2$ AGH University of Krakow,
\small 
al. Mickiewicza 30, 30-059 Krak\'ow, 
\small Poland\\
\small $^3$ TU Bergakademie Freiberg,
\small 09596 Freiberg,
\small Germany\\
}

\begin{document}

\date{}

\maketitle 

\begin{abstract}
\noindent
Given an integer $\Delta \ge 3$, let ${\cal G}_{\Delta }$ be the set of connected graphs $G\neq K_{\Delta +1}$   
with maximum degree $\Delta $  and, for $i=1,\cdots, \Delta $, let  $V_i(G)$ be the set of vertices of $G$ of degree $i$. \\
We prove that $\sum\limits_{i=1}^\Delta c_i|V_i(G)|$ is a lower bound on the independence number $\alpha(G)$ of $G\in {\cal G}_\Delta$, where  $c_\Delta=\frac{1}{\Delta}$ and  $ic_{i}=1-c_{i+1}$ for $i=1,\cdots,\Delta-1$.
Moreover, if $\varepsilon >0$ and $j\in \{1,\cdots, \Delta\}$, then the inequality $\alpha(G)\ge \varepsilon|V_j(G)|+\sum\limits_{i=1}^\Delta c_i|V_i(G)|$  does not hold for infinitely many graphs $G\in {\cal G}_\Delta$.
We also show that an independent  set $I\subset V(G)$ of $G\in {\cal G}_\Delta$ such that $|I|\ge  \sum\limits_{i=1}^\Delta c_i|V_i(G)|$ can be found in polynomial time.\\ 
Finally,  further lower bounds on $\alpha(G)$ in terms of degrees of $G$ are presented. 
\end{abstract}

\section{Introduction and Results}\label{I}

We consider simple, finite, and undirected graphs $G$, where  
 $V(G)$, $E(G)$, $N_G(v)$, and \\$d_G(v)=|N_G(v)|$ denote the vertex set of $G$, the edge set of $G$, the neighbourhood of $v\in V(G)$ in $G$, and the degree of $v\in V(G)$ in $G$, respectively.  For terminology and notation not defined here 
we refer to \cite{West}.\\
A subset $I$ of $V(G)$ is called {\it independent} if the subgraph of $G$ induced by $I$ is edgeless. 
The {\it independence number $\alpha(G)$} of $G$ is the  cardinality of a largest
independent set of $G$. \\~\\An independent set  of a graph $G$ with cardinality $ \alpha(G)$ is called a
{\it maximum independent set} of $G$. Finding a maximum independent set in a given graph $G$ is known to be  NP-hard. For this reason, it is interesting to determine  lower bounds on the independence number of a graph.\\
Lower bounds on $\alpha(G)$  in terms of degrees of $G$ are proved in \cite{ACL,C,H,HR,HS,KP,M,Wei}. 

For $\Delta \ge 3$, let ${\cal G}_{\Delta }$ be the set of connected graphs $G\neq K_{\Delta +1}$   
with maximum degree $\Delta $  and, for $i=1,\cdots, \Delta $, let  $V_i(G)$ be the set of vertices of $G$ of degree $i$.\\
An immediate consequence of the famous and well-known Theorem of Brooks (\cite{B}) on the relationship between the maximum degree of a graph and its chromatic number is the following inequality (\ref{E1}).\\ 
If $G\in {\cal G}_\Delta$, then
\begin{equation}\label{E1}
\alpha(G) \geq \frac{1}{\Delta}|V(G)|.
\end{equation}
Moreover, there are infinitely many $\Delta$-regular graphs $G\in {\cal G}_\Delta$ such that (\ref{E1}) is tight.

 Here we  ask the following question. \\Do there exist positive real numbers $c_1,\cdots,c_{\Delta}$ such that 
$\alpha(G)\ge \sum\limits_{i=1}^\Delta c_i|V_i(G)|$ is true for all  $G\in {\cal G}_\Delta$   and  if $\varepsilon >0$ and  $j\in \{1,\cdots,\Delta\}$, then the inequality 
$
\alpha(G)\ge \varepsilon|V_j(G)|+\sum\limits_{i=1}^\Delta c_i|V_i(G)|
$  does not hold  for infinitely many graphs $G\in {\cal G}_\Delta$ ?

Our paper is organized as follows. Our following results Theorem \ref{T1} in case $\Delta'=\Delta$ and Theorem \ref{T2} give a positive answer to the above question. Further lower bounds on $\alpha(G)$  are presented by  Theorem \ref{T1} in case $\Delta'<\Delta$ and by Theorem \ref{T3}. The proofs  are postponed to Section \ref{proofs}.

\begin{theorem}\label{T1}~~\\
If $3\le \Delta'\le \Delta $,   $c_\Delta=\frac{1}{\Delta}$,   $ic_{i}+c_{i+1}=1$ for $i=1,\cdots,\Delta-1$, and $G\in{\cal G}_{\Delta'}$, then
\begin{equation}\label{E2}
\alpha(G)\ge \sum\limits_{i=1}^{\Delta'}  c_i|V_i(G)|
\end{equation} 
and an independent set $I$ of $G$ such that $|I|\ge \sum\limits_{i=1}^{\Delta'} c_i|V_i(G)|$ can be found in polynomial time. 
\end{theorem}

For example, if  $\Delta'=4$ and $\Delta=6$, then (\ref{E2}) has the form\\ 
$\alpha(G)\ge \frac{91}{144}|V_1(G)|+\frac{53}{144}|V_2(G)|+\frac{19}{72}|V_3(G)|+\frac{5}{24}|V_4(G)|$, whereas, in case $\Delta'=\Delta=4$, Theorem \ref{T1} leads to 
$\alpha(G)\ge \frac{5}{8}|V_1(G)|+\frac{3}{8}|V_2(G)|+\frac{1}{4}|V_3(G)|+\frac{1}{4}|V_4(G)|$.\\
Note that  $\frac{91}{144}>\frac{5}{8}$ and  $\frac{53}{144}<\frac{3}{8}$, hence, these two lower bounds on $\alpha(G)$ are not comparable, even though the coefficients $\frac{5}{8},\frac{3}{8},\frac{1}{4}$, and $\frac{1}{4}$ are best possible in the sense of the following Theorem \ref{T2}.

\begin{theorem}\label{T2}~\\
If $\Delta \ge 3$, $c_\Delta=\frac{1}{\Delta}$,  $ic_{i}+c_{i+1}=1$ for $i=1,\cdots,\Delta-1$, $\varepsilon >0$, and  $j\in \{1,\cdots,\Delta\}$, then the inequality 
$
\alpha(G)\ge \varepsilon|V_j(G)|+\sum\limits_{i=1}^\Delta c_i|V_i(G)|
$  does not hold  for infinitely many graphs $G\in {\cal G}_\Delta$.
\end{theorem}

Even if $\Delta'$ is small, the calculation of the coefficients $c_1,\cdots,c_{\Delta'}$ in inequality (\ref{E2}) of Theorem \ref{T1} by  the recursion $c_\Delta=\frac{1}{\Delta}$ and  $ic_{i}=1-c_{i+1}$ for $i=1,\cdots,\Delta-1$ is time-consuming  if $\Delta$ is  large. This situation is not improving significantly if the complicated explicit formulas  (b) and  (c)  of the forthcoming Lemma \ref{L1} are used. This is in contrast to the calculation of the coefficients $d_1,\cdots,d_{\Delta'}$ in inequality (\ref{E4}) of the following Theorem \ref{T3} by the recursion $d_1=1-\frac{1}{e}$ and $d_{i+1}=1-id_i$ for $i=1,\cdots,\Delta'-1$. 

\begin{theorem}\label{T3}
Let $\Delta'\ge 3$, $G\in {\cal G}_{\Delta'}$, $i\in \{1,\cdots,\Delta'\}$, \\
$d_i=\frac{1}{i+1}+(i-1)!\big(\frac{1}{e}-\sum\limits_{j=0}^{i+1}\frac{(-1)^{j}}{j!}\big)$ if $i$ is even, and 
$d_i=\frac{1}{i+1}+(i-1)!\big(\sum\limits_{j=0}^{i+1}\frac{(-1)^{j}}{j!}-\frac{1}{e}\big)$ if $i$ is odd, where $e$ denotes the Euler number. Then
\begin{equation}\label{E4}
\alpha(G)\ge \sum\limits_{i=1}^{\Delta'}d_i|V_i(G)|
\end{equation}
and an independent set $I$ of $G$ such that $|I|\ge \sum\limits_{i=1}^{\Delta'} d_i|V_i(G)|$ can be found in polynomial time. 
Moreover, the coefficients $d_1,\cdots,d_{\Delta'}$ fulfil the recursion $d_1=1-\frac{1}{e}$ and $d_{i+1}=1-id_i$ for $i=1,\cdots,\Delta'-1$.
\end{theorem}

Note that 
$\frac{5}{8}<\frac{91}{144} <  d_1=1-\frac{1}{e}$ and 
 $\frac{1}{4}<\frac{19}{72}<d_3=1-\frac{2}{e}$ for the coefficients $d_1$ and $d_3$ in Theorem \ref{T3}.
  Hence, if $\Delta'=4$ and $|V_1(G)\cup V_3(G)|-|V_2(G)\cup V_4(G)|$ is large enough, then the lower bound on $\alpha(G)$ in Theorem \ref{T3} is stronger than the previous bounds \\$\frac{91}{144}|V_1(G)|+\frac{53}{144}|V_2(G)|+\frac{19}{72}|V_3(G)|+\frac{5}{24}|V_4(G)|$ and
$\frac{5}{8}|V_1(G)|+\frac{3}{8}|V_2(G)|+\frac{1}{4}|V_3(G)|+\frac{1}{4}|V_4(G)|$.

\section{Proofs}\label{proofs}

Given $G\in {\cal G}_\Delta$, let in the sequel   $f_G(v)=c_i$ if $v\in V_i(G)$ for $i=1,\cdots,\Delta$.\\
First we list some properties of the coefficients $c_1,\cdots,c_\Delta$ in  Theorem \ref{T1} in the forthcoming Lemma \ref{L1}.

\begin{lemma}\label{L1}~\\
If $\Delta \ge 3$, $c_\Delta=\frac{1}{\Delta}$, and  $ic_{i}+c_{i+1}=1$ for $i=1,\cdots,\Delta-1$, then

{\em (a)} $c_{\Delta-1}=\frac{1}{\Delta}$, $c_1<1$, $c_{i+1}< c_i$ for $i=1,\cdots,\Delta-2$,\\
{\em (b)} if $i\le \Delta-2$ and $i\equiv \Delta \mod 2$, then
$c_{i}=\frac{1}{i+1}+\big(\sum\limits_{j=1}^{\frac{\Delta-i-2}{2}}\frac{i+2j}{i(i+1)\cdots (i+2j+1)}\big)+\frac{1}{i(i+1)\cdots \Delta}$, and\\
{\em (c)} if $i\le \Delta-3$ and $i\equiv \Delta-1 \mod 2$, then $c_i=\frac{1}{i+1}+\big(\sum\limits_{j=1}^{\frac{\Delta-i-3}{2}}\frac{i+2j}{i(i+1)\cdots (i+2j+1)}\big)+\frac{\Delta-1}{i(i+1)\cdots \Delta}  $.
\end{lemma}

\begin{proof}[{\bf Proof of Lemma \ref{L1}.}]~\\
Proof of  (a).\\
Because $(\Delta-1)c_{\Delta-1}+c_\Delta=1$ it follows $c_{\Delta-1}=\frac{1}{\Delta}$.\\
It holds $c_i=\frac{1-c_{i+1}}{i}=\frac{1-\frac{1-c_{i+2}}{i+1}}{i}=\frac{1}{i+1}+\frac{c_{i+2}}{i(i+1)}$ for $i\le \Delta -2$. Using
 $c_\Delta=c_{\Delta -1}=\frac{1}{\Delta}$, we have $c_i>0$ for all  $i\le \Delta $. Thus, $c_1<c_1+c_2=1$. Moreover, $c_i>\frac{1}{i+1}$ for $i\le \Delta-2$. Now assume $c_{i+1}\ge c_i$ for some $i\in \{1,\cdots,\Delta-2\}$. It follows $c_i=\frac{1-c_{i+1}}{i}<\frac{1-c_{i}}{i}$ implying $c_i\le \frac{1}{i+1}$, a contradiction. \\ 
The proofs of (b) and (c) are by induction on $i$.\\
Proof of (b).\\
If $i=\Delta-2$, then $\frac{\Delta-i-2}{2}=0$, hence, $\sum\limits_{j=1}^{\frac{\Delta-i-2}{2}}\frac{i+2j}{i(i+1)\cdots (i+2j+1)}=0$ and it follows\\
$c_{\Delta-2}=\frac{1}{\Delta-1}+\frac{c_{\Delta}}{(\Delta-2)(\Delta-1)}=\frac{1}{\Delta-1}+\frac{1}{(\Delta-2)(\Delta-1)\Delta}$.

Now let $\Delta \ge 5$ and $i\le \Delta-4$.
We have $c_i=\frac{1}{i+1}+\frac{c_{i+2}}{i(i+1)}$ and, by induction,\\
$c_{i+2}=\frac{1}{i+3}+\big(\sum\limits_{j=1}^{\frac{\Delta-i-4}{2}}\frac{i+2+2j}{(i+2)(i+3)\cdots (i+2j+3)}\big)+\frac{1}{(i+2)(i+3)\cdots \Delta}$. Hence,\\ \\
$c_i=\frac{1}{i+1}+\frac{c_{i+2}}{i(i+1)}=\frac{1}{i+1}+\frac{1}{i(i+1)(i+3)}+\frac{1}{i(i+1)}\big(\sum\limits_{j=1}^{\frac{\Delta-i-4}{2}}\frac{i+2+2j}{(i+2)(i+3)\cdots (i+2j+3)}\big)+\frac{1}{i(i+1)\cdots \Delta}   $\\
$=\frac{1}{i+1}+\frac{i+2}{i(i+1)(i+2)(i+3)}+\frac{1}{i(i+1)}\big(\sum\limits_{j=2}^{\frac{\Delta-i-2}{2}}\frac{i+2j}{(i+2)(i+3)\cdots (i+2j+1)}\big)+\frac{1}{i(i+1)\cdots \Delta}  $\\
$=\frac{1}{i+1}+\big(\sum\limits_{j=1}^{\frac{\Delta-i-2}{2}}\frac{i+2j}{i(i+1)\cdots (i+2j+1)}\big)+\frac{1}{i(i+1)\cdots \Delta}  $.

Proof of (c).\\
If $i=\Delta-3$, then  $\frac{\Delta-i-3}{2}=0$, hence,  $\sum\limits_{j=1}^{\frac{\Delta-i-3}{2}}\frac{i+2j}{i(i+1)\cdots (i+2j+1)}=0$  and it follows\\
$c_{\Delta-3}=\frac{1}{\Delta-2}+\frac{c_{\Delta-1}}{(\Delta-3)(\Delta-2)}=\frac{1}{\Delta-2}+\frac{\Delta-1}{(\Delta-3)(\Delta-2)(\Delta-1)\Delta}$.

Now let $\Delta \ge 6$ and $i\le \Delta-5$.
We have $c_i=\frac{1}{i+1}+\frac{c_{i+2}}{i(i+1)}$ and, by induction,\\ \\
$c_{i+2}=\frac{1}{i+3}+\big(\sum\limits_{j=1}^{\frac{\Delta-i-5}{2}}\frac{i+2+2j}{(i+2)(i+3)\cdots (i+2j+3)}\big)+\frac{\Delta-1}{(i+2)(i+3)\cdots \Delta}$.
Hence,\\ \\
$c_i=\frac{1}{i+1}+\frac{c_{i+2}}{i(i+1)}=\frac{1}{i+1}+\frac{1}{i(i+1)(i+3)}+\frac{1}{i(i+1)}\big(\sum\limits_{j=1}^{\frac{\Delta-i-5}{2}}\frac{i+2+2j}{(i+2)(i+3)\cdots (i+2j+3)}\big)+\frac{\Delta-1}{i(i+1)\cdots \Delta}   $\\
$=\frac{1}{i+1}+\frac{i+2}{i(i+1)(i+2)(i+3)}+\frac{1}{i(i+1)}\big(\sum\limits_{j=2}^{\frac{\Delta-i-3}{2}}\frac{i+2j}{(i+2)(i+3)\cdots (i+2j+1)}\big)+\frac{\Delta-1}{i(i+1)\cdots \Delta}  $\\
$=\frac{1}{i+1}+\big(\sum\limits_{j=1}^{\frac{\Delta-i-3}{2}}\frac{i+2j}{i(i+1)\cdots (i+2j+1)}\big)+\frac{\Delta-1}{i(i+1)\cdots \Delta}  $ and Lemma \ref{L1} is proved.
\end{proof}

 The minimum degree of a graph $H$ is denoted by $\delta(G)$  and, if there is no risk of confusion, we use $\delta$ instead of $\delta(H)$.   
 
 \begin{lemma}\label{L2}~\\
Let $3\le \Delta'\le \Delta$,  $c_\Delta=\frac{1}{\Delta}$,  $ic_{i}+c_{i+1}=1$ for $i=1,\cdots,\Delta-1$, and $G\in {\cal G}_{\Delta'}$ with $\delta(G)\le \Delta'-1$.\\
If  $H$ is a non-empty induced subgraph of $G$ ($H$ not necessarily connected, possibly $H=G$) and $\delta=\delta(H)$, then 
$$V_\delta^*(H)=\{u \in V_\delta(H) ~|~d_G(u)+\sum\limits_{w\in N_H(u)}d_G(w)\ge \delta^2+\delta+1\} \neq \emptyset .$$
Moreover, if $u\in V_\delta^*(H)$, then 
 $f_G(u)+\sum\limits_{w\in N_H(u)}f_G(w)\le 1.$ 
\end{lemma}
\begin{proof}[{\bf Proof of Lemma \ref{L2}}]~\\
Clearly, $\delta \le \delta(G)\le \Delta'-1$.
\\
First let $\delta=0$. Then $N_H(u)=\emptyset$ and, because $G$ is connected, $d_G(u)\ge 1$  for $u\in V_\delta(H)$, hence, $V_\delta^*(H)= V_\delta(H) \neq \emptyset$ and $f_G(u)\le c_1<1$ for $u\in V_\delta(H)$ by Lemma \ref{L1} (a) in this case.\\
Now let  $\delta \ge 1$. 
Clearly,   $N_H(v)\neq \emptyset$ and 
$d_G(v)\ge d_H(v)=\delta$ for $v\in V_\delta(H)$. 
Moreover, if $v\in V_\delta(H)$, then  $d_G(w)\ge d_H(v)\ge \delta$ for $w\in N_H(v)$.\\
For contradiction, we now assume that $V_\delta^*(H)=\emptyset$. \\Then $(\delta+1)\delta \le d_H(u)+\sum\limits_{w\in N_H(u)}d_H(w)\le d_G(u)+\sum\limits_{w\in N_H(u)}d_G(w)\le \delta^2+\delta$ for all $u\in V_\delta(H)$. Thus, $d_H(u)=d_G(u)=\delta$ and $d_H(w)=d_G(w)=\delta$ if $w\in N_H(u)$ and  $u\in V_\delta(H)$, hence, $N_H(u)\subset V_\delta(H)$ if $u\in V_\delta(H)$. \\
If $C$ is a component of $H$ with $V_\delta(H)\cap V(C)\neq \emptyset$, then $V(C)\subseteq V_\delta(H)$ and $d_G(u)=\delta $ for all $u\in V(C)$.\\
Because $V_{\Delta'}(G) \neq \emptyset$ and $G$ is connected, it follows $G\neq C$ and that there is an edge $uv$ with $u\in V(C)$ and $v\in V(G)\setminus V(C)$. Thus, $d_G(u)\ge d_H(u)+1= \delta+1$ for that $u$. This is a contradiction and $V_\delta^*(H)\neq \emptyset$ is proved.\\
We have $u\in V_\delta^*(H)$ if and only if $u \in V_\delta(H)$ and  $d_G(u)+\sum\limits_{w\in N_H(u)}d_G(w)\ge \delta^2+\delta+1$. This inequality holds if   at least on vertex in  $\{u\}\cup N_H(u)$ has degree at least $\delta+1$ in $G$, hence, $f_G(u)+\sum\limits_{w\in N_H(u)}f_G(w)\le \delta c_\delta+c_{\delta+1}=1$ because $\{c_1,\cdots, c_\Delta \}$ is decreasing  by Lemma \ref{L1} (a)  and $c_\Delta=c_{\Delta-1}=\frac{1}{\Delta}$. 
\end{proof}

 The   algorithm in the following Lemma \ref{L3} is a slight modification of a well-known greedy procedure for  generating independent sets  in graphs (\cite{M}). 

 \begin{lemma}\label{L3}~\\
 Let $\Delta' \ge 3$ and $G\in {\cal G}_{\Delta'}$ such that $\delta(G)\le \Delta'-1$.\\
 Then the following algorithm generates an independent set $\{u_1,...,u_k\}$ of $G$ in polynomial time.\\
 ALGORITHM:\\
$G_1:=G$ and  $j:=1$.\\
While $V(G_j)\neq \emptyset$ choose $u_j \in V_{\delta(G_j)}^*(G_j)$,  $G_{j+1}=G_j-(\{u_j\}\cup N_{G_j}(u_j))$,  and $j:=j+1$.\\ 
$k:=j-1$.\\  STOP

\end{lemma}
\begin{proof}[{\bf Proof of Lemma \ref{L3}}]~\\
It holds $\delta (G_j)\le \Delta -1$ for $j=2,\cdots,k$ and,
by Lemma \ref{L2},  $V_{\delta(G_j)}^*(G_j)\neq \emptyset$ for $j=1,\cdots, k$. Obviously,   $\{u_1,...,u_k\}$ is an  independent set of $G$ and it is generated in polynomial time.
\end{proof}

\begin{proof}[{\em {\bf Proof of Theorem \ref{T1}}}]~\\
If $G$ is $\Delta'$-regular, then nothing needs to be proven because
there are fast algorithms to construct a $\Delta'$-colouring of $G$ (\cite{S}) and it holds $\alpha(G)\ge |I| \ge \frac{|V(G)|}{\Delta}$ for a largest colour class $I$ of such a colouring.\\ 
Hence, we may assume that  $\delta(G)\le \Delta'-1$.
Let $\{u_1,...,u_k\}$ be an independent set of $G$ obtained by the algorithm of Lemma \ref{L3}.
By Lemma \ref{L2}, 
 $f_G(u_j)+\sum\limits_{w\in N_{G_j}(u)}f_G(w)\le 1$ for $j=1,\cdots, k$ and it follows\\
$\sum\limits_{i=1}^{\Delta'} c_i|V_i(G)|=\sum\limits_{v\in V(G)}f_G(v)=\sum\limits_{j=1}^k \bigg(f_G(u_j)+\sum\limits_{w\in N_{G_j}(u)}f_G(w)\bigg)\le k$. 
\end{proof}

\begin{proof} [{\bf Proof of Theorem \ref{T2}}]~\\ Let $\varepsilon>0$ be given and we will show by constructions that Theorem \ref{T2} holds for all \\$j\in \{1,\cdots,\Delta\}$.

Let $G'$ be a connected  $\Delta$-regular graph with $|V(G')|=k$ and let the graph $G$ be constructed in two steps as follows. For $x\in V(G')$ let $H_x$ be a copy of $K_\Delta$ and let the graph $G$ be the disjoint union of these $H_x$ for all $x\in V(G')$, i.e. $|V(G)|=k\Delta$. Next, for every edge $xy\in E(G')$ add an edge $e(xy)$ connecting a vertex of $H_x$ with a vertex of $H_y$ to $E(G)$, such that in the resulting graph $G$ the edges $e(xy)$ for $xy\in E(G')$ form a matching. Then $G$ is also  
$\Delta$-regular and $G\neq K_{\Delta+1}$, hence, $G\in {\cal G}_\Delta$.
Clearly, $G$  contains $k$ mutually disjoint copies of $K_\Delta$, $V_\Delta(G)=V(G)$, and $\alpha(G)=k$ by the mentioned result of Brooks. Thus,
$\alpha(G)=k=\frac{1}{\Delta }|V_{\Delta}(G)|$ and if
 an infinite set of such graphs $G$ is considered, then it follows that Theorem \ref{T2} holds for $j= \Delta$.

Now let  the infinite set ${\cal A}_\Delta \subset {\cal G}_\Delta$ be defined in two steps as follows.\\
1. $K_\Delta \in {\cal A}_\Delta$\\
2. If $G'\in {\cal A}_\Delta$ and $x\in V_{\Delta -1}(G')$, then the graph $G$ obtained from the disjoint union of $G'$ and $K_\Delta$ and adding the edge $xy$, where $y\in V(K_\Delta)$, belongs to ${\cal A}_\Delta $.\\ 
If $G\in {\cal A}_\Delta$ and step 2 of the definition of ${\cal A}_\Delta$  is used $k-1$ times to construct $G$, then  $|V(G)|=k\Delta$, $V(G)=V_{\Delta -1}(G)\cup V_\Delta(G)$, $|V_{\Delta -1}(G)|=k\Delta-(2k-2)$, and $|V_{\Delta }(G)|=(2k-2)$. Since $G$ contains $k$ mutually disjoint copies of $K_\Delta$ and $G$ is $k$- colorable by Brooks' Theorem, it follows $\alpha(G)=k$.
Thus,
$\alpha(G)=k=c_{\Delta -1}|V_{\Delta -1}|+c_{\Delta }|V_{\Delta }|=\frac{1}{\Delta}|V_{\Delta -1}|+\frac{1}{\Delta}|V_{\Delta }|$ and we are done in case $j= \Delta-1$.

It remains  the case $j\le \Delta-2$. We consider 
 the infinite set ${\cal B}_j \subset {\cal G}_\Delta$  defined as follows:\\
If $G'\in {\cal A}_\Delta$, then for every $x\in V_{\Delta -1}(G')$ let a disjoint copy  of $H_x=K_{j+1}$   and the edge $xy$ be added to $G'$, where $y\in V(H_x)$. Then the resulting graph $G$ belongs to ${\cal B}_j$.\\
If $|V(G')|=k\Delta$, then, using $V(G')=V_{\Delta -1}(G')\cup V_\Delta(G')$, $|V_{\Delta -1}(G')|=k\Delta-(2k-2)$, and $|V_{\Delta }(G')|=(2k-2)$, it follows $V(G)=V_{j}(G)\cup V_{j+1}(G)\cup V_\Delta(G)$, \\$|V_{j}(G)|=(k\Delta-(2k-2))j$, $|V_{j+1}(G)|=k\Delta-(2k-2)$, and $|V_{\Delta }(G)|=k\Delta$. \\Clearly, $\alpha(G)=k+(k\Delta-(2k-2))$. \\Since $c_{\Delta}=\frac{1}{\Delta}$ we obtain
$\alpha(G)=c_{j}|V_{j}(G)|+c_{j+1 }|V_{j+1}(G)|+\frac{1}{\Delta }|V_{\Delta }(G)|$ because $jc_{j}+c_{j+1}=1$, thus, Theorem \ref{T2} is  proved.
\end{proof}

\begin{proof}[{\em {\bf Proof of Theorem \ref{T3}}}]~\\
Let $\Delta > \Delta'$ and  let the coefficients $c_1,\cdots,c_\Delta$ be defined as is Theorem \ref{T1}.\\
First we use the properties (a),  (b), and (c) of Lemma \ref{L1}. \\
Given $\varepsilon >0$, let $\Delta $ be chosen large enough such that 
$c_{i}> \frac{1}{i+1}+\sum\limits_{j=1}^{\infty}\frac{i+2j}{i(i+1)\cdots (i+2j+1)}-\varepsilon$ \\for $i=1,\cdots ,\Delta'$.\\
Next note that $\sum\limits_{j=1}^{\infty}\frac{i+2j}{i(i+1)\cdots (i+2j+1)}=\sum\limits_{j=1}^{\infty}\frac{1}{i(i+1)\cdots (i+2j)}-\sum\limits_{j=1}^{\infty}\frac{1}{i(i+1)\cdots (i+2j+1)}$.\\
If $i\in \{1,\cdots,\Delta'\}$  is even, then \\
$\sum\limits_{j=1}^{\infty}\frac{1}{i(i+1)\cdots (i+2j)}=(i-1)!\sum\limits_{j=1}^{\infty}\frac{1}{(i+2j)!}=(i-1)!\big(\sum\limits_{j=0}^{\infty}\frac{1}{(2j)!}-\sum\limits_{j=0}^{\frac{i}{2}}\frac{1}{(2j)!}\big)$\\
$=(i-1)!\big(\frac{1}{2}(e+\frac{1}{e})-\sum\limits_{j=0}^{\frac{i}{2}}\frac{1}{(2j)!}\big)$\\
and
$\sum\limits_{j=1}^{\infty}\frac{1}{i(i+1)\cdots (i+2j+1)}=(i-1)!\sum\limits_{j=1}^{\infty}\frac{1}{(i+2j+1)!}=(i-1)!\big(\sum\limits_{j=0}^{\infty}\frac{1}{(2j+1)!}-\sum\limits_{j=0}^{\frac{i}{2}}\frac{1}{(2j+1)!}\big)$\\$
=(i-1)!\big(\frac{1}{2}(e-\frac{1}{e})-\sum\limits_{j=0}^{\frac{i}{2}}\frac{1}{(2j+1)!}\big)$.

If $i\in \{1,\cdots,\Delta'\}$  is odd, then\\ 
$\sum\limits_{j=1}^{\infty}\frac{1}{i(i+1)\cdots (i+2j)}=(i-1)!\sum\limits_{j=1}^{\infty}\frac{1}{(i+2j)!}=(i-1)!\big(\sum\limits_{j=0}^{\infty}\frac{1}{(2j+1)!}-\sum\limits_{j=0}^{\frac{i-1}{2}}\frac{1}{(2j+1)!}\big)$\\
$=(i-1)!\big(\frac{1}{2}(e-\frac{1}{e})-\sum\limits_{j=0}^{\frac{i-1}{2}}\frac{1}{(2j+1)!}\big)$\\
and
$\sum\limits_{j=1}^{\infty}\frac{1}{i(i+1)\cdots (i+2j+1)}=(i-1)!\sum\limits_{j=1}^{\infty}\frac{1}{(i+2j+1)!}=(i-1)!\big(\sum\limits_{j=0}^{\infty}\frac{1}{(2j)!}-\sum\limits_{j=0}^{\frac{i+1}{2}}\frac{1}{(2j)!}\big)$\\
$=(i-1)!\big(\frac{1}{2}(e+\frac{1}{e})-\sum\limits_{j=0}^{\frac{i+1}{2}}\frac{1}{(2j)!}\big)$.

Hence, if $i\in \{1,\cdots,\Delta'\}$ is even, then  \\
$c_i+\varepsilon>\frac{1}{i+1}+ (i-1)!\big(\frac{1}{2}(e+\frac{1}{e})-\sum\limits_{j=0}^{\frac{i}{2}}\frac{1}{(2j)!}\big)-(i-1)!\big(\frac{1}{2}(e-\frac{1}{e})-\sum\limits_{j=0}^{\frac{i}{2}}\frac{1}{(2j+1)!}\big)$\\
$=\frac{1}{i+1}+(i-1)!\big(\frac{1}{e}-\sum\limits_{j=0}^{\frac{i}{2}}\frac{1}{(2j)!}+\sum\limits_{j=0}^{\frac{i}{2}}\frac{1}{(2j+1)!}\big)=\frac{1}{i+1}+(i-1)!\big(\frac{1}{e}-\sum\limits_{j=0}^{i+1}\frac{(-1)^{j}}{j!}\big)=d_i$ and,

if $i\in \{1,\cdots,\Delta'\}$ is odd, then \\
$c_i+\varepsilon>\frac{1}{i+1}+(i-1)!\big(\frac{1}{2}(e-\frac{1}{e})-\sum\limits_{j=0}^{\frac{i-1}{2}}\frac{1}{(2j+1)!}\big)-(i-1)!\big(\frac{1}{2}(e+\frac{1}{e})-\sum\limits_{j=0}^{\frac{i+1}{2}}\frac{1}{(2j)!}\big)$\\
$=\frac{1}{i+1}+(i-1)!\big(-\frac{1}{e}-\sum\limits_{j=0}^{\frac{i-1}{2}}\frac{1}{(2j+1)!}+\sum\limits_{j=0}^{\frac{i+1}{2}}\frac{1}{(2j)!}\big)=\frac{1}{i+1}+(i-1)!\big(\sum\limits_{j=0}^{i+1}\frac{(-1)^{j}}{j!}-\frac{1}{e}\big)=d_i$. 

By  Theorem \ref{T1}, 
$\alpha(G)\ge \sum\limits_{i=1}^{\Delta'} c_i|V_i(G)|> \sum\limits_{i=1}^{\Delta'} (d_i-\varepsilon)|V_i(G)|$ for all $\varepsilon >0$ implying \\$\sum\limits_{i=1}^{\Delta'} c_i|V_i(G)|\ge \sum\limits_{i=1}^{\Delta'} d_i|V_i(G)|$ and an independent set $I$ of $G$ such that $|I|\ge \sum\limits_{i=1}^{\Delta'} c_i|V_i(G)|$ can be found in polynomial time.\\
To complete the proof of Theorem \ref{T3}, it is checked easily that $d_1=1-\frac{1}{e}$ and  $d_{i+1}=1-id_i$ for $i=1,\cdots,\Delta'-1$.
\end{proof}

 \section*{Declarations}
No data have been used. There are no competing interests. All authors contributed equally.

\end{document}